\documentclass{gtart_a}
\pdfoutput=1
\usepackage{pinlabel}

\title{Distortion in transformation groups}

\author{Danny Calegari}
\givenname{Danny}
\surname{Calegari}
\address{Department of Mathematics\\
California Institute of Technology\\
Pasadena CA 91125\\
USA}
\email{dannyc@its.caltech.edu}
\urladdr{}

\author{Michael H Freedman}
\givenname{Michael H}
\surname{Freedman}
\address{Microsoft Research\\
1 Microsoft Way\\
Redmond WA 98052\\
USA}
\email{michaelf@microsoft.com}
\urladdr{}

\def\appendixauthor#1{\author{Appendix by #1}}
\appendixauthor{Yves de Cornulier}
\givenname{Yves}
\surname{de Cornulier}
\address{{\rm Yves de Cornulier:}\ \ Institut de Math\'ematiques\\
Universit\'e de Neuch\^atel\\\newline
Rue \'Emile Argand 11\\
CH-2007 Neuch\^atel\\
Switzerland}
\email{yves.cornulier@ens.fr}
\urladdr{}

\volumenumber{10}
\issuenumber{}
\publicationyear{2006}
\papernumber{7}
\lognumber{0669}
\startpage{267}
\endpage{293}

\doi{}
\MR{}
\Zbl{}

\arxivreference{math.DS/0509701}
\arxivpassword{j3pbg}


\received{7 October 2005}
\revised{20 February 2006}
\accepted{8 February 2006}
\published{26 March 2006}
\publishedonline{26 March 2006}
\proposed{Benson Farb}
\seconded{Leonid Polterovich, Robion Kirby}
\corresponding{}
\editor{}
\version{}

\let\xysavmatrix\xymatrix
\def\xymatrix{\disablesubscriptcorrection\xysavmatrix}
\AtBeginDocument{\let\bar\wbar}
\def\theshortauthors{Danny Calegari and Michael H Freedman}

\makeatletter
\def\cnewtheorem#1[#2]#3{\newtheorem{#1}{#3}[section]
\expandafter\let\csname c@#1\endcsname\c@thm}
\newtheorem{thm}{Theorem}[section]
\cnewtheorem{lem}[thm]{Lemma}
\cnewtheorem{cor}[thm]{Corollary}
\cnewtheorem{conj}[thm]{Conjecture}
\cnewtheorem{qn}[thm]{Question}
\cnewtheorem{prob}[thm]{Problem}
\cnewtheorem{mprob}[thm]{Meta-Problem}

\newtheorem{Athm}{Theorem}

\makeautorefname{Athm}{Theorem}

\newtheorem*{S2_thm}{Theorem A}
\newtheorem*{S1_thm}{Theorem B}
\newtheorem*{Sn_thm}{Theorem C}

\theoremstyle{definition}
\cnewtheorem{defn}[thm]{Definition}
\cnewtheorem{construct}[thm]{Construction}
\cnewtheorem{note}[thm]{Notation}
\cnewtheorem{conv}[thm]{Convention}
\makeautorefname{defn}{Definition}

\theoremstyle{remark}
\cnewtheorem{rmk}[thm]{Remark}
\cnewtheorem{exa}[thm]{Example}
\makeatother

\def\R{\mathbb R}
\def\N{\mathbb N}

\def\C{\mathbb C}
\def\Z{\mathbb Z}
\def\Q{\mathbb Q}

\def\SS{\mathbf S}
\def\GL{\textnormal{GL}}

\def\homeo{\textnormal{Homeo}}
\def\diff{\textnormal{Diff}}

\def\supp{\textnormal{supp}}

\def\MCG{\textnormal{MCG}}
\def\LAT{\textnormal{LAT}}
\def\LONG{\textnormal{LONG}}
\def\id{\textnormal{Id}}
\def\til{\wwtilde}
\def\Gal{\textnormal{Gal}}

\begin{document}

\begin{asciiabstract}
We exhibit rigid rotations of spheres as distortion elements in groups
of diffeomorphisms, thereby answering a question of J Franks and M
Handel.
We also show that every homeomorphism of a sphere
is, in a suitable sense, as distorted as possible
in the group Homeo(S^n), thought of as a discrete group.

An appendix by Y de Cornulier shows that Homeo(S^n) has the
strong boundedness property, recently introduced by G
Bergman. This means that every action of the discrete group
Homeo(S^n) on a metric space by isometries has bounded orbits.
\end{asciiabstract}

\begin{htmlabstract}
<p class="noindent">
We exhibit rigid rotations of spheres as distortion elements in groups of
diffeomorphisms, thereby answering a question of J Franks and M Handel.
We also show that every homeomorphism of a sphere is, in a suitable sense,
as distorted as possible in the group Homeo(<b>S</b><sup>n</sup>),
thought of as a discrete group.
</p>
<p class="noindent">
An appendix by Y de Cornulier shows that Homeo(<b>S</b><sup>n</sup>)
has the <em>strong boundedness</em> property, recently introduced
by G Bergman. This means that every action of the discrete group
Homeo(<b>S</b><sup>n</sup>) on a metric space by isometries has
bounded orbits.
</p>
\end{htmlabstract}

\begin{abstract}
We exhibit rigid rotations of spheres as distortion elements in groups of
diffeomorphisms, thereby answering a question of J Franks and M Handel.
We also show that every homeomorphism of a sphere is, in a suitable sense,
as distorted as possible in the group $\mathrm{Homeo}(\mathbf{S}^n)$,
thought of as a discrete group.

An appendix by Y de Cornulier shows that $\mathrm{Homeo}(\mathbf{S}^n)$
has the {\em strong boundedness} property, recently introduced
by G Bergman. This means that every action of the discrete group
$\mathrm{Homeo}(\mathbf{S}^n)$ on a metric space by isometries has
bounded orbits.
\end{abstract}

\maketitle

\section{Introduction}
\label{sec1}

The study of abstract groups as geometric objects has a long history, but has
been pursued especially vigorously since the work of Gromov
\cite{Gromov_hyperbolic,Gromov_AIDG}.  Typically the focus
is on finitely presented groups; however, interesting results have also
been obtained from this perspective in the theory of
transformation groups --- ie groups of homeomorphisms
of manifolds.

The topic of this paper is distortion in transformation groups, especially
groups of homeomorphisms of spheres.
Informally, an element $h$ in a finitely generated group $G$ is {\em distorted}
if the word length of $h^n$ grows sublinearly in $n$. One also sometimes says that
the {\em translation length} of $h$ vanishes.
Geometrically, this corresponds to the condition that the homomorphism from
$\Z$ to $G$ sending $n$ to $h^n$ is not a quasi-isometric embedding.

One can also make sense of the concept of distortion in infinitely generated groups.
An element $h$ in a (not necessarily finitely generated) group $G$ is distorted
if there is a finitely generated subgroup $H$ of $G$ containing $h$
such that $h$ is distorted in $H$ as above.
To show that an element is undistorted, one typically tries to define an appropriate
real-valued function on $G$ which is (almost) subadditive, and which grows
linearly on $h^n$. For example, {\em quasi-morphisms} are useful in this respect, and
highlight one point of contact between distortion and the theory of bounded cohomology.
On the other hand, exhibiting distortion is typically done {\it ad hoc}, and there do
not seem to be many very general or flexible constructions known.

In this paper, we study distortion in groups of homeomorphisms of spheres,
especially groups consisting of transformations with a definite amount
of analytic regularity (ie $C^1$ or $C^\infty$). By contrast with
Polterovich \cite{Polterovich}, Franks--Handel \cite{FranksHandel} or
Ghys--Gambaudo \cite{GhysGambaudo}, we do not insist that our groups
preserve a probability measure; the considerable additional flexibility
this affords has the consequence that our results have more of an
existential character than those of the papers cited above, exhibiting
distortion rather than ruling it out.

\subsection{Statement of results}

\begin{note}
The letters $G,H$ will denote groups of some sort, and $S$ a (symmetric) generating set,
although $\SS^n$ denotes the $n$--sphere.
If $G$ is a group, and $H$ is a subgroup, we write $H < G$. The group
$G$ will often be a transformation group on some manifold, and a typical
element $h$ will be a homeomorphism or diffeomorphism of some analytic quality.
The letters $i,j,n$ will denote integers, $r$ will denote a degree of smoothness,
and $g$ will denote a growth function (ie a function $g\co\N \to \N$).
$c$ and $k$ will usually denote (implicit) constants in some inequality.
We let $\R_+$ denote the {\em non-negative} real numbers.
Other notation will be introduced as needed.
\end{note}

In \fullref{sec2} and \fullref{sec3} we summarize some basic definitions and study examples of
distorted and undistorted elements in various groups.

In \fullref{sec4} we exhibit rigid rotations of $\SS^2$ as distortion elements in the
group of $C^\infty$ diffeomorphisms of the sphere.

Our main result in this section is:

\begin{Athm}
\label{thma}
For any angle $\theta \in [0,2\pi)$ the rigid rotation $R_\theta$ of $\SS^2$
is a distortion element in a finitely generated subgroup of $\diff^\infty(\SS^2)$.
Moreover, the distortion function of $R_\theta$ can be chosen to grow
faster than any given function.
\end{Athm}

Here $R_\theta$ is a clockwise rotation about a fixed axis through angle
$\theta$. To say that the distortion function grows faster than any given
function means that for any $g\co\N \to \N$ we can find a finitely generated
group $G<\diff^\infty(\SS^2)$
for which there are words of length $\sim n_i$ in the generating
set which express powers $R_\theta^{f(n_i)}$ of $R_\theta$ for some sequence
$n_i \to \infty$, where $f(n)>g(n)$ for all sufficiently large $n \in \N$.
In this case we say that the distortion function grows {\em faster than $g$}.

This answers a question of John Franks and Michael Handel, motivated by results in their
paper \cite{FranksHandel}.

In \fullref{sec5} we go down a dimension, and study rigid rotations of $\SS^1$. Our main
result here is:

\begin{Athm}
\label{thmb}
For any angle $\theta \in [0,2\pi)$ the rigid rotation $R_\theta$ of $\SS^1$ is
a distortion element in a finitely generated subgroup of $\diff^1(\SS^1)$. 
Moreover, the distortion function of $R_\theta$ can be chosen to grow faster
than any given function.
\end{Athm}

The proof of \fullref{thmb} makes use of Pixton's results from \cite{Pixton}, and
the arguments should be familiar to people working in the theory of
foliations. It should be remarked that our construction cannot be made $C^2$,
and it appears to be unknown whether a rigid rotation of 
$\SS^1$ is distorted in $\diff^\infty(\SS^1)$ (or even in $\diff^2(\SS^1)$).

\begin{rmk}
The possibility of proving \fullref{thmb} was pointed out to the first author by Franks and
Handel, after reading an early version of this paper.
\end{rmk}

In \fullref{sec6} we relax our analytic conditions completely, and study distortion
in the full group of homeomorphisms of $\SS^n$. Here our main result is quite general:

\begin{Athm}
\label{thmc}
Fix $n\geq 1$.
Let $h_1, h_2, \dots$ be any countable subset of $\homeo(\SS^n)$, and
$g_1,g_2,\dots \co \N \to \N$ any countable collection of growth functions. Then
there is a finitely generated subgroup $H$ of $\homeo(\SS^n)$ (depending on $\lbrace h_i \rbrace$
and $\lbrace g_i \rbrace$)
such that {\em every} $h_i$ is simultaneously distorted in $H$. Moreover, the
distortion function of $h_i$ grows faster than $g_i$.
\end{Athm}

The proof of \fullref{thmc} uses the full power of the Kirby--Siebenmann theory of
homeomorphisms of manifolds for a key factorization step.
It is an interesting question whether one can
exhibit distortion in an arbitrary homeomorphism of the sphere without recourse
to such sophisticated technology.

Finally, in an appendix, Yves de Cornulier uses the proof of
\fullref{thmc} to show that the group $\homeo(\SS^n)$ is {\em strongly
bounded}. Here an abstract group $G$ is said to be strongly
bounded if every symmetric subadditive non-negative real-valued
function on $G$ is bounded. A countable group has this property if
and only if it is finite.

\subsection{Acknowledgements}

The first author would like to thank Michael Handel for
suggesting the problem which motivated \fullref{thma}, and to thank him
and John Franks for reading preliminary versions of this paper,
and for making clarifications and corrections. He would also
like to thank Daniel Allcock for some useful comments.

\section{Distortion elements}
\label{sec2}

\subsection{Conjugation notation}

\begin{note}
For a group $G$ and elements $a,b \in G$, we abbreviate
the conjugate $b^{-1}ab$ by
$$a^b:=b^{-1}ab.$$
\end{note}

Notice with this convention that
$$(a^b)^c = a^{bc}.$$

\subsection{Basic definitions}

\begin{defn}
Let $G$ be a finitely generated group, and let $S$ be a finite generating set.
By convention, we assume $S=S^{-1}$. Given $h \in G$, the {\em length of $h$
with respect to $S$} is the minimum integer $n$ such that $h$ can be
expressed as a product
$$h = s_1s_2 \cdots s_n$$
where each $s_i \in S$. We write
$$\ell_S(h)=n.$$
By convention, we take $\ell_S(1)=0$.
\end{defn}

Note that $\ell_S$ is a {\em subadditive function}; that is, for all $h_1,h_2 \in G$,
$$\ell_S(h_1h_2) \le \ell_S(h_1) + \ell_S(h_2).$$
Moreover, it is {\em non-negative} and {\em symmetric}; ie $\ell(h) = \ell(h^{-1})$.
This motivates the definition of a {\em length function} on a group $G$.

\begin{defn}
Let $G$ be a group. A {\em length function} on $G$ is a function $L\co G \to \R_+$ satisfying
$L(1)=0$ which is symmetric and subadditive.
\end{defn}

The function $\ell_S$ depends on the choice of generating set $S$, but only up to
a multiplicative constant:

\begin{lem}\label{rewriting_bound}
If $S_1,S_2$ are two finite generating sets for $G$, then there is a constant
$c \ge 1$ such that
$$\frac 1 c \ell_{S_2}(h) \le \ell_{S_1}(h) \le c\ell_{S_2}(h)$$
for all $h \in G$.
\end{lem}
\begin{proof}
Each $s \in S_1$ can be expressed as a word of length $n(s)$ in the elements of $S_2$, and
vice versa. Then take $c$ to be the maximum of the $n(s)$ over all $s \in S_1 \cup S_2$.
\end{proof}

\begin{defn}\label{distortion_definition}
Let $G$ be a finitely generated group, and let $S$ be a symmetric finite generating
set as above. The {\em translation length} of an element
$h \in G$, denoted $\|h\|_S$, is the limit
$$\|h\|_S:=\lim_{n \to \infty} \frac {\ell_S(h^n)} n.$$
An element $h \in G$ is a {\em distortion element} if the translation length is $0$.
\end{defn}

\begin{rmk}
Note that by the subadditivity property of $\ell_S$, the limit exists. Moreover,
by \fullref{rewriting_bound}, the property of being a distortion element is
independent of the choice of generating set $S$.
\end{rmk}

\begin{rmk}
With this definition, torsion elements are distortion elements. Some authors
(including Franks and Handel \cite{FranksHandel}) explicitly require distortion
elements to be {\em nontorsion}.
\end{rmk}

Sometimes, we shall pay attention to the growth rate of $\ell_S(h^n)$ as a function
of $n$ to make qualitative distinctions between different kinds of distortion
elements. If $h$ is not torsion, we define the {\em distortion function} to be the function
$$D_{S,h}\co \N \to \N$$
defined by the property
$$D_{S,h}(n) = \max \lbrace i \; | \; \ell_S(h^i) \le n \rbrace.$$
We can remove the dependence of this function on $S$ as follows.
For two functions
$$f,g\co \N \to \N$$
we write $f \precsim g$ if there is a constant $k\ge 1$ such that
$$f(n) \le kg(kn+k) + k \; \text{ for all } n \in \N$$
and then write $f \sim g$ if $f \precsim g$ and $g \precsim f$.
It is straightforward to see that $\precsim$ is transitive, and that
$\sim$ is an equivalence relation.
In case $f \sim g$, we say that $f,g$ are {\em quasi-equivalent}.
With this definition, the quasi-equivalence class of $D_{S,h}$ is independent of
$S$, and may be denoted $D_h$.

We are also interested in comparing growth rates in a cofinal sense:
\begin{defn}
Given $g\co \N \to \N$ we say that the distortion function of $h \in G$ 
(with respect to a generating set $S$) {\em grows faster than $g$} if 
there is a sequence $n_i \to \infty$ and a function $f\co \N \to \N$ such
that $f(n) > g(n)$ for all sufficiently large $n$, and such that
$$\ell_S(h^{f(n_i)}) \le n_i.$$
\end{defn}

We say for example that $h$ has {\em quadratic distortion} if $g(n)=n^2$ or
{\em exponential distortion} if $g(n) = e^n$ as above.

Finally, we may define a distortion element in an arbitrary group:

\begin{defn}
Let $G$ be a group. An element $h \in G$ is a distortion element if there is
a finitely generated subgroup $H < G$ with $h \in H$ such that $h$ is a distortion
element in $H$.
\end{defn}

Note that for such an element $h$, the quasi-equivalence class of the
distortion function may certainly depend on $H$.

\subsection{Examples}

\begin{exa}
In $\Z$ only the identity element is distorted.
\end{exa}

\begin{exa}
If $\phi\co G \to H$ is a homomorphism, and $\phi(h)$ is not distorted in $H$, then $h$
is not distorted in $G$.
\end{exa}

\begin{exa}
If $L\co G \to \R_+$ is a length function, and 
$$\lim_{n \to \infty} \frac {L(h^n)} n > 0$$ then $h$ is not
distorted. More generally, a length function gives a lower bound for
word length with respect to any finite generating set, and therefore
an upper bound on distortion. For example, if $L(h^n)$ grows like
$\log(n)$ then $h$ is no more than exponentially distorted.
\end{exa}

The next few examples treat distortion in linear groups.

\begin{exa}\label{firstGL}
Let $G = \GL(n,\C)$ and define $L\co G \to \R_+$ by
$$L(A) = \log \text{ of the $\max$ of the operator norms of $A$ and $A^{-1}$}.$$
Then $L$ is a length function. It follows that if $A$ has an eigenvalue
with absolute value $\ne 1$ then $A$ is not distorted.
\end{exa}

\begin{exa}\label{secondGL}
Let $\sigma \in \Gal(\C/\Q)$ be a Galois automorphism of $\C$. Then $A$ is distorted in $\GL(n,\C)$
if and only if $\sigma(A)$ is. It follows that if $A$ is distorted, then every eigenvalue
must be algebraic, with all conjugates on the unit circle.
\end{exa}

\begin{exa}
Let $G < \GL(n,\C)$ be a finitely generated subgroup with entries in a number field $K$.
We may construct length functions from valuations associated to finite primes in the
ring of integers of $K$. If $x \in K$ then $v(x) = 0$ for all discrete valuations $v$ on $K$
if and only if $x$ is a unit. A unit in a number field with absolute value $1$ is a root
of unity; cf Lang \cite{Lang}. Combined with \fullref{firstGL} and \fullref{secondGL},
one can show that an arbitrary element $A \in \GL(n,\C)$ is distorted if and only
if every eigenvalue of $A$ is a root of unity. Note that the distortion of a non-torsion
element is at most {\em exponential}. See Lubotzky--Mozes--Raghunathan \cite{LubotzkyMozesRaghunathan} for details.
\end{exa}

\begin{exa}
In the Baumslag--Solitar group $\langle a,b \; | \; aba^{-1} = b^2 \rangle$ the element
$b$ has exponential distortion. Similarly, in the group
$$\langle a,b,c \; | \; aba^{-1} = b^2, bcb^{-1} = c^2 \rangle$$
the element $c$ has doubly-exponential distortion. Note that as a corollary, we deduce
that this second group is not linear. This example and others are
mentioned by Gromov \cite[Chapter 3]{Gromov_AIDG}.
\end{exa}

\begin{exa}
Let $G$ be a group. A {\em quasi-morphism} is a map
$\phi\co G \to \R$ such that there is a constant $c>0$ for which
$$|\phi(h_1) + \phi(h_2) - \phi(h_1h_2)| \le c$$
for all $h_1,h_2 \in G$. If $|\phi(h)|>c$ then $h$ is not distorted.
\end{exa}

Quasi-morphisms are intimately related to (second) bounded cohomology.
See eg Ghys--Gambaudo \cite{GhysGambaudo} for a salient discussion.

\section{Distortion in transformation groups}
\label{sec3}

\subsection{Transformation groups}

\begin{note}
For a compact $C^\infty$ manifold $M$, we denote the group of homeomorphisms of $M$ by $\homeo(M)$,
and the group of $C^r$ diffeomorphisms by $\diff^r(M)$, where $r = \infty$ is possible.
Here a homeomorphism $h$ is in $\diff^r(M)$ if both it and its inverse are $C^r$. Note
that this implies $dh$ has full rank everywhere. If we wish
to restrict to orientation-preserving subgroups, we denote this by a $+$
superscript.
\end{note}

\subsection{Distortion in $\diff^1$}\label{general_nonsense}

Suppose $M$ is a smooth compact Riemannian manifold, and $h \in \diff^1(M)$.
We define the following norm:
$$\|dh\| : = \log \sup_{v \in UTM} |dh(v)|$$
where $|dh(v)|$ denotes the length of $dh(v)$, and the supremum is taken over
all vectors $v$ in the unit tangent bundle of $M$.

Note that since $h$ is a diffeomorphism and $M$ is compact, $dh$ cannot be strictly
contracting at every point, and therefore $\|dh\|\ge 0$. If we define
$$\|dh\|^+ = \max (\|dh\|, \|d(h^{-1})\|)$$ 
then it is clear that $\|d\cdot\|^+$ is a length function on $\diff^1(M)$.
In general, therefore, the growth rate of $\|dh^n\|^+$ as a function of $n$
puts an upper bound on the distortion function of $h$ in any finitely generated subgroup
of $\diff^1(M)$.

Note if we choose two distinct Riemannian metrics on $M$, the length functions $\|d\cdot\|^+$ they
define will be quasi-equivalent, by compactness. On the other hand, if $M$ is non-compact,
different quasi-isometry classes of Riemannian metrics may give rise to 
qualitatively different length functions.

\begin{exa}
Suppose $h$ has a fixed point $p$ and $dh|_{TpM}$ has an eigenvalue with
absolute value $\ne 1$. Then $h$ is not distorted in $\diff^1(M)$.
\end{exa}

\begin{exa}\label{positive_Lyapunov}
Oseledec's theorem (see Pollicott \cite[Chapter 2]{Pollicott})
says that for $h \in \diff^1(M)$ where $M$ is
a compact manifold, and for $\mu$ an ergodic $h$--invariant
probability measure on $M$, there are real numbers $\lambda_1 > \cdots > \lambda_k$
called {\em Lyapunov exponents},
and a $\mu$--measurable $dh$--invariant splitting $TM = \oplus_{i=1}^k E^i$
such that
$$\lim_{n \to \infty} \frac 1 n \log |dh^n(v)| = \lambda_l$$
for almost every $v \in \oplus_{i=l}^k E^i$ but not in $\oplus_{i=l+1}^k E^i$.
In particular, if $\lambda_1>0$, then $h$ is undistorted in $\diff^1(M)$.
\end{exa}

\begin{exa}\label{positive_entropy}
Let $M$ be a compact manifold, and suppose $h \in \diff^1(M)$ has positive
topological entropy. Then there is an ergodic $h$--invariant probability measure
$\mu$ for which $h$ has positive $\mu$--entropy.
The Pesin--Ruelle inequality (see \cite[Chapter 3]{Pollicott}) says
$$\sum_{\lambda_i>0} \lambda_i \ge \mu\text{-entropy of } h$$
where the $\lambda_i$ are the Lyapunov exponents for $h$ with respect to the
measure $\mu$. It follows that some Lyapunov exponent $\lambda_1$ for
$\mu$ is positive, and therefore, as in \fullref{positive_Lyapunov},
$h$ is undistorted in $\diff^1(M)$.
\end{exa}

By contrast, if $\|dh^n\|^+$ is bounded independently of $n$, 
then the group $\langle h \rangle$ is equicontinuous,
and is precompact in the group of Lipschitz homeomorphisms of $M$, by the
Arzela--Ascoli theorem.
By Repov\u{s}--\u{S}\u{c}epin \cite{RS} (ie the Hilbert--Smith conjecture for Lipschitz actions),
a compact group of Lipschitz homeomorphisms of a smooth manifold $M$
is a Lie group. In our case, this group is abelian, since it contains a dense
abelian subgroup $\langle h \rangle$, and is therefore (up to finite index)
a finite dimensional torus. Thus the uniformly equicontinuous case reduces to that of
torus actions.

A key case to understand in this context is when the torus in question is $\SS^1$, and
the simplest example is that of a rigid rotation of a sphere. It is this
example which we study in the next few sections.

\section{Rotations of $\SS^2$}\label{S2_section}
\label{sec4}

\subsection{The group $G$}\label{group_G}

We describe a particular explicit group $G < \diff^\infty(\SS^2)$ which will
be important in the sequel. By stereographic projection, we may identify
$\SS^2$ conformally with $\C \cup \infty$.

Let $T$ be the similarity
$$T\co z \to 2z.$$
Then $\langle T \rangle$ acts discretely and properly discontinuously on $\C^*$ with quotient
a (topological) torus. A fundamental domain for the action is the annulus $A$ defined by
$$A = \lbrace z \in \C \; | \; 1 \le |z| \le 2 \rbrace.$$
We let $\partial A^+$ and $\partial A^-$ denote the components $|z| = 2$ and
$|z|=1$ of $\partial A$ respectively. We define a disk $D$ by
$$D = \lbrace z \in \C \; | \; |z-3/2| \le 1/4 \rbrace.$$
We let $F$ be a $C^\infty$ diffeomorphism with the following properties:
\begin{itemize}
\item{$F$ is the identity outside the annulus $0.99 \le |z| \le 2.01$}
\item{$F$ restricted to the annulus $1.01 \le |z| \le 1.99$ agrees with the rotation $z \to -z$}
\end{itemize}

We define $G = \langle T,F \rangle$, and think of it as a subgroup
of $\diff^\infty(\SS^2)$ fixing $0$ and $\infty$. Notice that for every $h \in G$ either
$h(D)$ is disjoint from $D$, or else $h|_D = \id|_D$. If $G^D$ denotes the stabilizer of $D$
in $G$, then we may identify the orbit $GD$ with the product $D \times S$ where $S$
is the set of (right) cosets of the subgroup $G^D$ in $G$. Note that $S$ is a
set with a (left) $G$--action. This action determines the action of $G$ on $D \times S$.

An explicit set of coset representatives for $S$ is the set of elements of the form $T^n$ and
$FT^n$ for all $n \in \Z$.

\subsection{Wreath products}\label{wreath_subsection}

Let $G,S$ and $D \subset \SS^2$ be as in \fullref{group_G}.
Let $\zeta_t, t \in \R$ be a 1--parameter subgroup of diffeomorphisms
of the unit disk with support contained in the interior. After conjugating
by a diffeomorphism, we think of $\zeta_t$ as a 1--parameter subgroup of
$\diff^\infty(\SS^2)$ with support contained in the interior of $D$.

\begin{defn}
Let $\R^S$ denote the set of functions from $S$ to $\R$, which can be thought
of as an abelian group with respect to addition.
The {\em wreath product} $G \wr_S \R$ is the semi-direct product
$$0 \to \R^S \to G \wr_S \R \to G \to 0$$
where $G$ acts on $\R^S$ by
$$f^h(s) = f(hs)$$
for $h \in G, s\in S$.
\end{defn}

The choice of 1--parameter group $\zeta_t$ determines a faithful homomorphism
$$\rho\co G\wr_S \R \to \homeo(\SS^2)$$
as follows. For $f \in \R^S$, define
$$\rho(f) = \prod_{s \in S} \zeta_{f(s)}^{\overline{s}}$$
where $\overline{s} \in G$ is a coset representative of $s \in S$.
Together with the action of $G$ on $\SS^2$
(in its capacity as a transformation group) this defines a faithful homomorphism $\rho$.
For the sake of brevity, in the sequel we will omit $\rho$, and think of $G\wr_S\R$ itself
as a subgroup of $\homeo(\SS^2)$.

\subsection{Analytic quality}
Given $f \in \R^S$, thought of as an element of $\homeo(\SS^2)$
as in \fullref{wreath_subsection}, the analytic quality of $f$
is {\it a priori} only $C^0$. However, if we can estimate the $C^r$ norm
of $f(T^n),f(FT^n)$ as $|n| \to \infty$, we can improve this {\it a priori}
estimate.

Notice that any $f \in \R^S$ is $C^\infty$ away from $0,\infty$.
In particular, any $f$ with finite support is $C^\infty$ on all of $\SS^2$.
Furthermore, conjugation by $T$ preserves the $C^1$ norm, and blows up the
$C^r$ norm by $2^{r-1}$, whereas conjugation by $F$ preserves the $C^r$ norm for every $r$. 
It follows that if we have an estimate
$$|f(T^n)|,|f(FT^n)| = o(2^{-|n|(r-1)})$$
as $|n| \to \infty$, then $f$ is $C^r$ at $0$
(here our notation $|f(s)|$ just means the absolute value of $f(s)$ for $s \in S$,
where we think of $f$ as a function from $S$ to $\R$).
By the change of co-ordinates $z \to 1/z$ one sees that $f$ is also $C^r$
at $\infty$ under the same hypothesis, and is therefore $C^r$ on all of $\SS^2$.

We summarize this as a lemma:

\begin{lem}\label{analytic_quality}
Let $f \in \R^S$ be thought of as an element of $\homeo(\SS^2)$ as in
\fullref{wreath_subsection}. Then we have the following estimates:
\begin{itemize}
\item{If $|f(s)|$ is bounded independently of $s \in S$ then $f$ is Lipschitz}
\item{If $\lim_{s \to \infty} |f(s)| = 0$ then $f$ is $C^1$}
\item{If $|f(T^n)|,|f(FT^n)| \to 0$ faster than any exponential (as a function of $n$),
then $f$ is $C^\infty$}
\end{itemize}
\end{lem}

\subsection{Rotations of $\SS^2$}

For each $\theta \in [0,2\pi)$ we let $R_\theta$ denote the rigid rotation
of $\SS^2$ with fixed points equal to $0$ and $\infty$. In stereographic co-ordinates,
$$R_\theta\co z \to e^{i\theta}z$$
where $z \in \C \cup \infty$.
Notice that $R_\pi$ is just multiplication by $-1$.

For $\theta \in \pi\Q$ the element $R_\theta$ is torsion in $\diff^\infty(\SS^2)$.
We will show in this section that $R_\theta$ is a distortion element in
$\diff^\infty(\SS^2)$ for arbitrary $\theta$. Moreover, the distortion function can
be taken to grow faster than any given function.

\subsection{Factorizing rotations}\label{factorize_subsection}

We can factorize $R_\theta$ in a natural way
as a product of two diffeomorphisms whose support is contained in
closed subdisks of $\SS^2$. This will be important for some later applications.

Let $B$ (for bump) be a smooth function $B\co \R^+ \to [0,1]$ which satisfies the
following properties:

\begin{itemize}
\item{$B(t)=0$ for $t<1/2$ and $B(t)=1$ for $t>2$}
\item{$B(t) + B(1/t) = 1$}
\item{$B$ is monotone decreasing and strictly positive on $(1/2,2)$}
\item{$B$ is infinitely tangent to the constant function $1$ at $2$ and
to the constant function $0$ at $1/2$}
\end{itemize}

For $\theta \in \R$, define $R_\theta^+$ by
$$R_\theta^+\co z \to e^{i B(|z|)\theta}z$$
and define $R_\theta^-$ by the identity
$$R_\theta^+ R_\theta^- = R_\theta.$$
Notice that as $\theta$ varies over $\R$, the set of transformations $R_\theta^-$ and
$R_\theta^+$ form smooth subgroups of $\diff^+(\SS^2)$. Moreover, the support
of the group $\lbrace R_\theta^- \: | \: \theta \in \R \rbrace$ is equal to the
disk
$$E^- = \lbrace z \; | \; |z|\le2 \rbrace$$
Similarly, the support of $R_\theta^+$ is the disk (in $\SS^2$)
$$E^+ = \lbrace z \; | \: |z| \ge 1/2 \rbrace.$$
Notice the important fact that $z \to 1/z$ conjugates $R_{-\theta}^+$ to
$R_\theta^-$ for any $\theta$. The reason for the sign change is that
a 1--parameter family of rotations which has a clockwise sense at one fixed
point has an anticlockwise sense at the other fixed point.

\subsection{Construction of the group}\label{construction_subsection}

Throughout the remainder of this section we assume that $\theta$ has been
fixed.

We define a diffeomorphism $Z$ which takes care of some bookkeeping for us.
Basically, the diffeomorphism $Z$ lets us move back and forth between
the 1--parameter groups $R^\pm_t$ with support in $E^\pm$ and a 1--parameter group
$\zeta_t$ with support in $D$, as in \fullref{group_G} and \fullref{wreath_subsection}.
The exact details of how this is done are irrelevant, but we must make an
explicit choice, which accounts for the (annoying) notational complexity below.

Let $Z \in \diff^\infty(\SS^2)$ satisfy the following properties:
\begin{itemize}
\item{$Z$ takes $D$ to $E^-$ and conjugates $R_t^-$ to a 1--parameter
subgroup $\zeta_t$:
$$\zeta_t := (R_t^-)^Z$$}
\item{$Z$ takes $FD$ to $T^3E^+$ (ie the image of the disk $E^+$ under the
similarity $z \to 8z$) and conjugates $(R_t^+)^{T^{-3}}$ to $\zeta_{-t}^F$:
$$\zeta_{-t}^F = (R_t^+)^{T^{-3}Z}$$}
\end{itemize}

The existence of such a diffeomorphism $Z$ follows from the disjointness
of the disks $E^-,T^3E^+$ and the fact that the subgroups $R_t^-$ and
$R_{-t}^+$ are abstractly conjugate, by $z \to 1/z$, as pointed out in
\fullref{factorize_subsection}.

Now form the group $\R^S$ as in \fullref{wreath_subsection} by means of
the subgroup $\zeta_t = (R_t^-)^Z$.

Let $t_i \in \R$ be chosen for all non-negative integers $i$ subject to the following
constraints:
\begin{itemize}
\item{$t_i = n_i \theta \mod 2\pi$ where $n_i \to \infty$ grow as fast as desired
(ie faster than some growth function $g\co \N \to \N$ we are given in advance)}
\item{$t_i \to 0$ faster than any exponential function}
\end{itemize}

Define the element
$f \in \R^S$ by
$$f(T^i) = t_i \text{ if } i \ge 0, \; f(T^i) = 0 \text{ if } i< 0, \;
f(FT^i) = 0 \text{ for all } i.$$
By \fullref{analytic_quality} the function $f$ is in $\diff^\infty(\SS^2)$ with
respect to the identification of $\R^S$ with a subgroup of $\homeo(\SS^2)$.

Now, for any $i$, the element
$$f_i := f^{T^i} (f^{T^i F})^{-1}$$
is contained in $\R^S$, and satisfies
$$f_i(s) = \begin{cases}
t_i & \text{ if } s = \id \\
-t_i & \text{ if } s = F \\
0 & \text{ otherwise.}\\
\end{cases}$$
We conjugate the $f_i$ back by $Z^{-1}$, and define
$$h_i := f_i^{Z^{-1}}.$$
Then $h_i$ agrees with $R_{t_i}^-$ on $E^-$ and agrees with $(R_{t_i}^+)^{T^{-3}}$
on $T^3E^+$.

Notice that $h_i$ preserves the foliation of $\SS^2\backslash \lbrace 0,\infty \rbrace$
by circles of equal latitude, and acts on each of these circles by a rotation.
Let $\LAT< \diff^\infty(\SS^2)$ denote the group of diffeomorphisms with this property;
ie informally, $\LAT$ preserves latitude, and acts as a rotation on each circle
with a fixed latitude. An element of $\LAT$ can be thought of as a $C^\infty$ function
$$\text{ latitude } \to \text{ rotation angle }$$
up to constant functions with values in $2\pi\Z$,
and any element of $\LAT$ can be recovered pictorially from the graph of this function.
Notice that $h_i \in \LAT$. In this way, we can abbreviate $h_i$ by a picture:

\begin{figure}[ht]
\labellist\small
\pinlabel {rotation angle} [l] at 0 103
\pinlabel {latitude} [bl] at 580 0
\endlabellist
\centerline{\includegraphics[width=8cm]{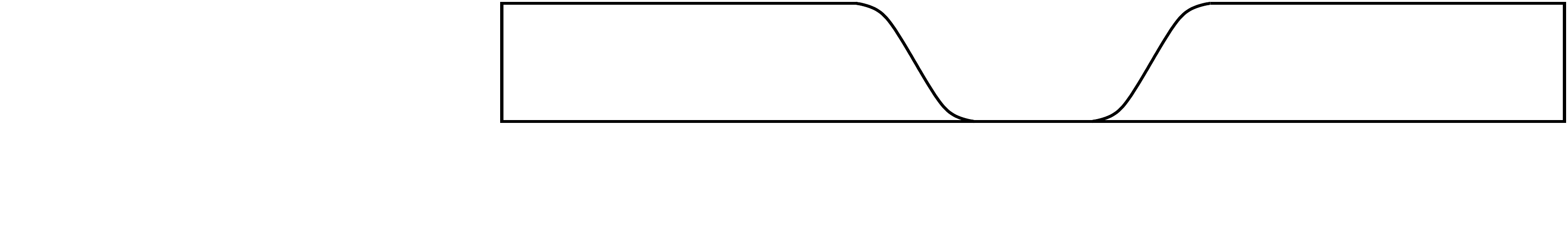}}
\caption{The element $h_i \in \LAT$ represented pictorially by the graph
  of a function}
\label{graph_figure}
\end{figure}

Let $\LONG < \diff^\infty(\SS^2)$ denote the group of diffeomorphisms of the form
$$z \to z\cdot u(|z|)$$
where $u\co \R^+ \to \R^+$ is infinitely tangent to the identity at $0$ and at $\infty$.
Informally, $\LONG$ is the group of diffeomorphisms which reparameterizes the set of latitudes, without
changing longitudes. Then $\LONG$ is contained in the normalizer of the group $\LAT$.
The conjugation action of $\LONG$ on $\LAT$ is given pictorially by reparameterizing
the base of the graph.

We claim that there are elements $M_1,M_2,M_3 \in \LONG$
such that for any $h_i$ we have an identity
$$h_i(h_i)^{M_1}(h_i)^{M_2}((h_i)^{M_3})^{-1} = R_{2t_i}.$$
The proof is given graphically by \fullref{conjugacy_figure}:

\begin{figure}[ht]
\labellist\small
\pinlabel {$+$} at 692 290
\pinlabel {$+$} at 1412 290
\pinlabel {$-$} at 170 65
\pinlabel {$=$} at 940 65
\endlabellist
\centerline{\includegraphics[width=12cm]{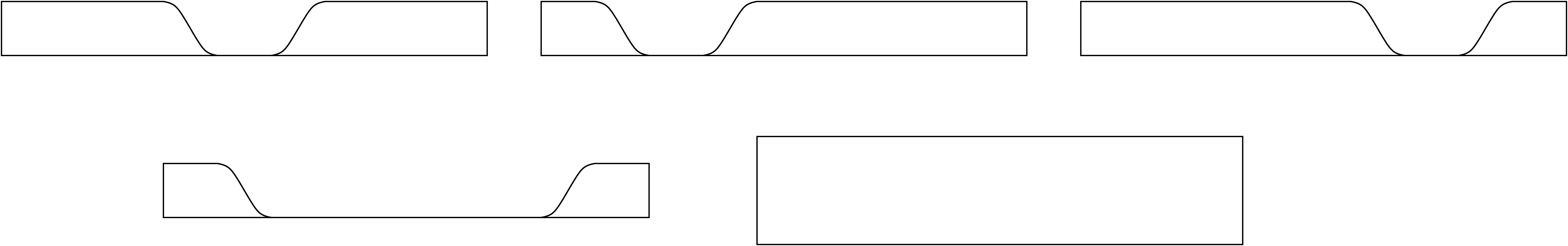}}
\caption{These figures denote the conjugates of $h_i$, and
demonstrate how an appropriate algebraic product of these conjugates is equal to
$R_{2t_i}$}\label{conjugacy_figure}
\end{figure}
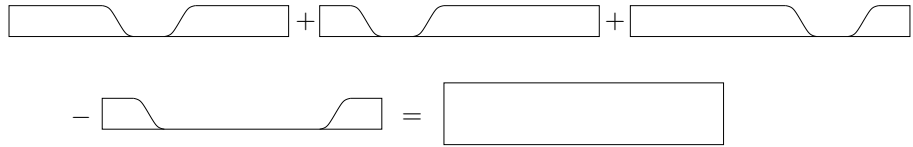

Now, $R_{2t_i} = R_\theta^{2n_i}$. Since the $n_i$ have been chosen to
grow faster than any function given in advance,
we have proved the following theorem:

\begin{S2_thm}
For any angle $\theta \in [0,2\pi)$ the rigid rotation $R_\theta$ of $\SS^2$
is a distortion element in a finitely generated subgroup of $\diff^\infty(\SS^2)$.
Moreover, the distortion function of $R_\theta$ can be chosen to grow
faster than any given function.
\end{S2_thm}

\section{Rotations of $\SS^1$}
\label{sec5}

In this section we show how to modify the construction of \fullref{S2_section}
to exhibit a rigid rotation as a distortion element in the group $\diff^1(\SS^1)$.
But first, we exhibit a rotation as a distortion element
in the group of {\em Lipschitz} homeomorphisms of $\SS^1$.

\subsection{Rotations of $\SS^1$}\label{S1_subsection}

As in the previous section, we denote by
$R_\theta$ the rotation of $\SS^1$ through angle $\theta \in [0,2\pi)$.

The first difference with \fullref{S2_section}
is that we cannot factorize a 1--parameter group of rotations
as the product of two 1--parameter groups with support contained in an interval.
(One way to see this is to use Poincar\'e's rotation number;
see eg Sinai \cite{Sinai} for a definition and basic properties.)

Let $\theta$ be fixed, and we choose $t_i \to 0,n_i \to \infty$ as $i \in \Z$ goes from
$0$ to $\infty$, with
$$t_i = n_i\theta \mod 2\pi$$
as in \fullref{construction_subsection}.

Let $I^\pm$ be two intervals which form an open cover of $\SS^1$. Then for
$t_i$ sufficiently close to $0$, we can factorize $R_{t_i}$ as a product of
two diffeomorphisms $\xi_i,\zeta_i$ with support contained in $I^+,I^-$
respectively. It is clear that we may choose $\xi_i,\zeta_i$ so that their
support is exactly equal to an interval, and they are both conjugate to
translations on these intervals.

Let $J$ be an open interval in $\SS^1$ which we parameterize by arclength
as $[-1,1]$. We let $T$ be a diffeomorphism of $\SS^1$ with support equal to
$J$, and with no fixed points in $J$. Then the restriction of $T$ to $J$
is conjugate to a translation, and we let $J_i$ for $i \in \Z$ be a tiling
of $J$ by fundamental domains for $t$.

Fix one such interval $J_0 \subset J$ and let $F$ be a diffeomorphism
of $\SS^1$ with support equal to $J_0$, and with no fixed points in $J_0$.
We let $J_{0i}$ for $i \in \Z$ be a tiling of $J_0$ by fundamental domains for $F$.

The group $G=\langle T,F \rangle$ acts as before on the set of translates of $J_{00}$,
and for all $h \in G$, either $h(J_{00})$ is disjoint from $J_{00}$, or
else $h|_{J_{00}} = \id|_{J_{00}}$. The interval $J_{00}$ is the analogue of
the disk $D$ from \fullref{S2_section}, and the elements $T,F$ are the analogues of
the diffeomorphisms of the same names in that section. The difference is that
if $GJ_{00}$ denotes the set of translates of $J_{00}$ by $G$, then $F|_{GJ_{00}}$
has infinite order, rather than order $2$.

Now let $Z^\pm$ be diffeomorphisms of $\SS^1$ taking $I^\pm$ respectively to
the interval $J_{00}$. The diffeomorphisms $Z^\pm$ are the analogue of the diffeomorphism
$Z$ from \fullref{S2_section}; the reason we need two such diffeomorphisms rather than
just one is that the factorization of $R_{t_i}$ into $\xi_i \zeta_i$ is no longer
canonical.

We let $f^+ \in \homeo(\SS^1)$ have support contained in $J$, and define it to be
the product
$$f^+ = \prod_{i = 0}^\infty \prod_{j = 0}^\infty \xi_i^{(Z^+)^{-1}F^{-j}T^{-i} }$$
and similarly, define
$$f^- = \prod_{i = 0}^\infty \prod_{j = 0}^\infty
\zeta_i^{(Z^-)^{-1}F^{-j}T^{-i}}.$$
Notice by \fullref{analytic_quality} that $f^\pm$ are Lipschitz (though {\em not}
$C^1$).

Then for each $i$,
$$\left( (f^+)^{T^i}((f^+)^{T^iF^{-1}})^{-1} \right)^{Z^+} = \xi_i$$
and
$$\left( (f^-)^{T^i}((f^+)^{T^iF^{-1}})^{-1} \right)^{Z^-} = \zeta_i$$
and therefore $R_{n_i\theta}$ can be expressed as a word of length $\sim 8i$
in the group
$$\langle f^+,f^-,F,T,Z^+,Z^-\rangle.$$
Notice that there is no analogue of the groups $\LAT$ and $\LONG$, and
consequently no analogue of the elements $M_1,M_2,M_3$.

\subsection{A $C^1$ example}\label{C1_subsection}

By a slight modification, using a trick of Pixton
we can actually improve the Lipschitz example of
\fullref{S1_subsection} to a $C^1$ example.

We note that by suitable choice of factorization of $R_{t_i}$ we can assume the
following:

\begin{itemize}
\item{The support of $\xi_i,\zeta_i$ is {\em exactly} equal to $I^+,I^-$ respectively}
\item{On $I^+$, each $\xi_i$ is conjugate to a translation, and similarly for
$I^-$}
\end{itemize}
Now, the elements $\xi_i$ for distinct $i$ will not be contained in a fixed
1--parameter subgroup of $\diff^\infty(I^+)$, but they are all {\em conjugate}
into a fixed 1--parameter subgroup, and similarly for the $\zeta_i$. 
The final condition we insist on is:
\begin{itemize}
\item{The conjugating maps can be taken to be
$C^1$ and converge in the $C^1$ topology to the identity.}
\end{itemize}

To see that this is possible, observe that for two $C^\infty$
diffeomorphisms $\epsilon$--close to the identity
in the $C^1$ norm, the commutator is $\epsilon^2$--close to the identity,
also in the $C^1$ norm. So for diffeomorphisms
$\phi_g$ defined by the property
$$\phi_g\co \theta \to \theta + g(\theta)$$
for $g\co \SS^1 \to \R$, we have that
$$\phi_{g_1}\phi_{g_2} \sim \phi_{g_1 + g_2}$$
with error which is comparable in size in the $C^1$ norm to the products of the
$C^1$ norms of $g_1,g_2$. Using this fact, one can readily produce a suitable
factorization.

\subsection{Pixton actions}\label{Pixton_action}

Consider an interval $I$ on which a diffeomorphism
$Y\co I \to I$ acts in a manner smoothly conjugate
to a translation, with fundamental domains $I_i$. Given another diffeomorphism
$\phi\co I_0 \to I_0$ we form the suspension $\Phi\co I \to I$ by
$$\Phi = \prod_i \phi^{Y^i}.$$
Note that $\langle Y,\Phi \rangle \cong \Z \oplus \Z$.
If $\phi$ restricted to $I_0$ is smoothly conjugate to a translation, then
{\it a priori} the action of $\langle Y,\Phi \rangle$ on $I$ is
Lipschitz. However, Pixton showed that it is {\em topologically} conjugate
(ie by a homeomorphism) to a $C^1$ action.

For the convenience of the reader, we give an outline of the construction of
a Pixton action. One
chooses co-ordinates on $I$ so that the ratio $|I_i|/|I_{i+1}|$ converges
to $1$ as $|i| \to \infty$. For instance, near $I$, the endpoints of the $I_n$
could be the harmonic series $1/2,1/3,\dots$ so that the ratio of successive lengths
is $i/(i+1) \to 1$. Then we require $Y\co I_i \to I_{i+1}$ to expand the linear
structure near the endpoints and contract it in the middle, so that the norm of
the first derivative of $\Phi|_{I_{i+1}}$ is smaller than that of
$\Phi|_{I_i}$ by a definite amount. Then both $Y$ and $\Phi$ are $C^1$ tangent to
the identity at the endpoints of $I$, and are therefore $C^1$ on the entire interval.
See Pixton \cite{Pixton} for rigorous details of this construction. One should remark that
a lemma of Kopell \cite{Kopell} implies that one cannot make the action $C^2$.

This construction has the following virtue:
if $\phi$ is contained in a smooth 1--parameter subgroup $\phi_t$, and
we form the associated 1--parameter subgroup $\Phi_t$ so that
$$\langle Y,\Phi_t \rangle \cong \Z \oplus \R$$
then we can form a Pixton action of this larger group which is $C^1$.

\subsection{$T$ and $X$}
Naively, one sees that by careful choice of $F$, one can arrange for the action of
$f^+|_{J_0}$ to be $C^1$. However, to make $f^+$ $C^1$ on all of $J$ requires us to
modify the definition slightly.

We will construct $X$, a diffeomorphism of $\SS^1$ with support equal to $J$,
conjugate to a translation on $J$, and with fundamental domains $J_i$, just
like $T$.

We let $\chi_t$ be a 1--parameter subgroup containing $\xi_0$. For each
$i$, we require that
$$\xi_i^{(Z^+)^{-1}T^{-i}X^i} \in (\chi_t)^{(Z^+)^{-1}}$$
which is possible, by the discussion at the end of
\fullref{C1_subsection}.
By choosing co-ordinates on $J$ suitably as above,
we can insist that both $X$ and $T$ are $C^1$.

Now we choose co-ordinates on $J_0$ so that $F$ and $(\chi_t)^{(Z^+)^{-1}}$ form a Pixton
action of $\Z \oplus \R$ there, as in \fullref{Pixton_action}.

We define
$$f^+ = \prod_{i=0}^\infty \prod_{j=0}^\infty
\xi_i^{(Z^+)^{-1}T^{-i}X^iF^{-j}X^{-i}}.$$
Note that $f^+$ is actually $C^1$.

Moreover, we have the following formula
$$\left( \left( (f^+)^{X^i}  ((f^+)^{X^iF^{-1}})^{-1} \right)^{X^{-i}T^i}
\right)^{Z^+} = \xi_i.$$
Relabelling $X$ as $X^+$ and defining $X^-$ similarly in terms of the $\zeta_i$,
one can define $f^-$ analogously. Putting this together, we have shown:

\begin{S1_thm}
For any angle $\theta \in [0,2\pi)$ the rigid rotation $R_\theta$ of $\SS^1$ is
a distortion element in a finitely generated subgroup of $\diff^1(\SS^1)$. 
Moreover, the distortion function of $R_\theta$ can be chosen to grow faster
than any given function.
\end{S1_thm}

One should remark that for a rigid rotation $R$ of $\SS^n$ where $n$ is arbitrary,
either $R$ has fixed points, in which case the construction of
\fullref{S2_section} shows that $R$ is a distortion element in $\diff^\infty(\SS^n)$,
or else the construction of this section can be generalized to show
that $R$ is a distortion element in $\diff^1(\SS^n)$,
in either case with distortion growing faster than any given function.

\begin{rmk}
Tsuboi showed in \cite{Tsuboi} that one can construct Pixton actions which
are $C^{1+\alpha}$ for every $\alpha < 1$. It is therefore likely that the
construction above exhibits a rigid rotation as an arbitrarily badly
distorted element in $\diff^{1+\alpha}(\SS^1)$.
\end{rmk}

By our discussion in \fullref{general_nonsense}, we make the following
conjecture:

\begin{conj}
Let $M$ be a compact smooth manifold, and let $h \in \diff^1(M)$. Then $h$ is a distortion element
in $\diff^1(M)$ whose distortion function can be chosen to grow faster than
any given function if and only if some finite power of $h$ is contained in
a Lipschitz action of a finite dimensional torus on $M$.
\end{conj}

Note that the ``only if" direction follows from \fullref{general_nonsense}.

\section{Distortion in $\homeo(\SS^n)$}
\label{sec6}

The group $\homeo(M)$ for an arbitrary manifold $M$ is considerably more complicated
than $\diff^\infty(M)$ or even $\diff^1(M)$. In this section, we first make a couple
of comments about distortion in $\homeo(M)$ in general, and then specialize to the
case of $\homeo(\SS^n)$.

\subsection{Mapping class groups}
For an arbitrary compact manifold $M$, there is a natural homomorphism
$$\homeo(M) \to \homeo(M)/\homeo_0(M) =: \MCG(M)$$
where $\homeo_0(M)$ is the normal subgroup consisting of homeomorphisms
isotopic to the identity, and $\MCG(M)$ is the {\em mapping class group}
of $M$. For reasonable $M$, this group is finitely presented, and quite
amenable to computation. Clearly for $h \in \homeo(M)$ to be a distortion
element, it is necessary for the image $[h]$ of $h$
in $\MCG(M)$ to be a distortion element.

\begin{exa}
A pseudo-Anosov homeomorphism of a closed surface $\Sigma$ of genus $\ge 2$ is not a
distortion element in $\homeo(\Sigma)$.
\end{exa}

\subsection{Distortion in $\homeo_0(M)$}

For suitable manifolds $M$, it is easy to find undistorted elements in $\homeo_0(M)$.

\begin{exa}\label{torus_wiggle}
Let $T^2$ denote the $2$--torus. Let $h\co T^2 \to T^2$ preserve the foliation of $T^2$
by meridians, and act as a rigid rotation on each meridian, where the angle of
rotation is not constant. This angle of rotation defines a map $\theta\co \SS^1 \to \SS^1$,
where the first factor labels the meridian, and the second factor is the amount of
rotation. If $\theta$ is homotopically trivial, $h$ is in $\homeo_0(T^2)$. In this
case, we claim $h$ is undistorted in $\homeo_0(T^2)$. To see this, suppose to the
contrary that $h$ is distorted in some finitely generated subgroup $H$.
Without loss of generality, we may expand $H$ to a larger finitely generated group,
where each generator $h_i$ has support contained in a closed disk in $T^2$.
If $\til{h_i}$ denotes a lift of $h_i$ to the universal cover $\R^2$, then there is
a constant $c$ such that
$$|d_{\R^2}(\til{h_i}(p),\til{h_i}(q)) - d_{\R^2}(p,q)| \le c$$
for any $p,q \in \R^2$. Without loss of generality, we may assume that the same
constant $c$ works for all $i$.

Now, if $I$ is a small transversal to the foliation of $T^2$ by meridians, intersecting
meridians where the function $\theta$ is nonconstant, it follows that
if we denote $I_n:= h^n(I)$, then a lift $\til{I_n}$ of $I_n$
has the property that the endpoints are
distance $\sim kn$ apart for some positive constant $k$.
By the discussion above, this implies that
any expression of $h^n$ in the generators $h_i$ and their inverses
has word length at least $\sim nk/c$. This
shows that $h$ is undistorted, as claimed.
\end{exa}

\begin{exa}
Let $M$ be a closed hyperbolic $3$--manifold. Let $\gamma \subset M$ be
a simple closed geodesic, and let $N$ be an embedded tubular neighborhood.
Let $h\co M \to M$ rotate $\gamma$ some distance, and be fixed outside
$N$. Then the argument of \fullref{torus_wiggle} shows that $h$ is
undistorted in $\homeo_0(M)$. Since $M$ is hyperbolic of dimension at
least $3$, Mostow rigidity \cite{Mostow}
implies that $\MCG(M)$ is finite. It follows that
$h$ is undistorted in the full group $\homeo(M)$.
\end{exa}

\begin{qn}
Is $h$ as in \fullref{torus_wiggle} undistorted in $\homeo(T^2)$?
\end{qn}

The method of construction in \fullref{torus_wiggle} produces an
undistorted element of $\homeo_0(M)$ when
$\pi_1(M)$ contains an undistorted element. Moreover, if $\MCG(M)$ is
finite, the element is undistorted in $\homeo(M)$. This begs the following obvious
question:

\begin{qn}\label{every_distorted}
Is there an infinite, finitely presented group $G$ in which every element is distorted?
\end{qn}

\begin{rmk}
A finitely presented infinite torsion group would answer \fullref{every_distorted}
affirmatively.
\end{rmk}

\begin{rmk}
It is worth observing that Ol'shanskii \cite{Olshanskii}
has shown the existence of a torsion-free
finitely generated group in which all elements are distorted, thereby answering
a question of Gromov.
\end{rmk}

The following construction gets around \fullref{every_distorted}, at a mild cost.

\begin{exa}
Let $M$ be a closed manifold with $\pi_1(M)$ infinite. Then $\til{M}$ inherits
a path metric pulled back from $M$ with respect to which the diameter is infinite.
It follows that $\til{M}$ contains a {\em ray} $r$ --- that is, an isometrically embedded
copy of $\R^+$ which realizes the minimal distance between any two points which it
contains. The ray $r$ projects to $M$ where it might intersect itself. By abuse
of notation, we refer to the projection as $r$. If the
dimension of $M$ is at least $3$, then we can perturb $r$ an arbitrarily small
amount so that it is embedded in $M$ (though of course not properly
embedded). In fact, we can even ensure that there is an embedded tubular neighborhood
$N$ of $r$ whose width tapers off to zero as one escapes to infinity in $r$
in its intrinsic path metric. Let $h$ be a homeomorphism of $M$, fixed outside $N$,
which translates the core (ie $r$) by some function
$$r(t) \to r(t+f(t))$$
where $f(t)$ is positive, and goes to $0$ as $t \to \infty$. Such a homeomorphism
may be constructed for instance by coning this translation of $r$ out to $\partial N$
with respect to some radial co-ordinates. Then $h$ might be distorted, but the
distortion function can be taken to increase as slowly as desired, by making
$f$ go to $0$ as slowly as desired. For example, we could ensure that the
distortion function grows slower than $n^\alpha$ for all $\alpha>1$.
\end{exa}

\subsection{Homeomorphisms of spheres}

We now specialize to $\SS^n$. We make use of the following seemingly innocuous lemma:

\begin{lem}[Kirby--Siebenmann, Quinn]\label{KS_lemma}
Let $h \in \homeo^+(\SS^n)$. Then $h$ can be factorized as a product
$$h = h_1 h_2$$
where the support of $h_1$ avoids the south pole, and the support of $h_2$ avoids the
north pole.
\end{lem}

For $h$ sufficiently close to the identity in the compact-open topology, this can
be proved by the geometric torus trick. For an arbitrary homeomorphism, it requires
the full power of topological surgery theory. See Kirby--Siebenmann \cite{KirbySiebenmann} for details
in the case $n\ne 4$ and Quinn \cite{Quinn} for the case $n=4$.

Using this lemma, we can produce another factorization:

\begin{lem}\label{factor_as_six}
Let $E_1,E_2$ be two closed disks in $\SS^n$ whose interiors cover $\SS^n$. Then
any $h \in \homeo^+(\SS^n)$ can be factorized as a product of at most 6 homeomorphisms,
each of which has support contained in either $E_1$ or $E_2$.
\end{lem}
\begin{proof}
Without loss of generality, we can assume that $E_1$ and $E_2$ contain collar
neighborhoods of the northern and southern hemisphere respectively.

Given $h \in \homeo^+(\SS^n)$, we factorize $h$ as $h_1 h_2$ as in \fullref{KS_lemma}.
Let $e_2$ be a radial expansion centered at the south pole, with
support contained in $E_2$, which takes
$\supp(h_1) \cap E_2$ into $E_2 \cap E_1$. Then $e_2 h_1 e_2^{-1}$ has support
contained in $E_1$. Similarly, we can find $e_1$ with support contained in $E_1$ such
that $e_1 h_2 e_1^{-1}$ has support contained in $E_2$.
Then
$$h = e_2^{-1} (e_2 h_1 e_2^{-1}) e_2 e_1^{-1} (e_1 h_2 e_1^{-1}) e_1$$
expresses $h$ as the product of 6 homeomorphisms, each with support in either $E_1$
or $E_2$.
\end{proof}

\begin{rmk}
Notice in the factorization in \fullref{factor_as_six} that the homeomorphisms
$e_1,e_2$ definitely depend on $h$.
\end{rmk}

\begin{Sn_thm}
Fix $n\geq 1$.
Let $h_1, h_2, \dots$ be any countable subset of $\homeo(\SS^n)$, and
$g_1,g_2,\dots \co \N \to \N$ any countable collection of growth functions. Then
there is a finitely generated subgroup $H$ of $\homeo(\SS^n)$ (depending on $\lbrace h_i \rbrace$
and $\lbrace g_i \rbrace$)
such that {\em every} $h_i$ is simultaneously distorted in $H$. Moreover, the
distortion function of $h_i$ grows faster than $g_i$.
\end{Sn_thm}
\begin{proof}
The subgroup $\homeo^+(\SS^n)$ of $\homeo(\SS^n)$ has index $2$, so after replacing
each $h_i$ by $h_i^2$ if necessary, we can assume each $h_i \in \homeo^+(\SS^n)$.

Fix a cover of $\SS^n$ by disks $E_1,E_2$ as in \fullref{factor_as_six}.
Let $n_i \to \infty$ grow sufficiently quickly, and relabel
the sequence
$$h_1^{n_1}, h_1^{n_2}, h_2^{n_2}, h_1^{n_3},h_2^{n_3},h_3^{n_3},h_1^{n_4},
\dots,h_1^{n_i},\dots,h_i^{n_i},h_1^{n_{i+1}}\dots$$
as $H_0,H_1,H_2,\dots$.

Applying \fullref{factor_as_six}, we write each
$H_i$ as a product
$$H_i = H_{i,1} H_{i,2} \dots H_{i,6}$$
where each $H_{i,j}$ has support contained in either $E_1$ or $E_2$.

From now on, the construction proceeds as in \fullref{sec4} and
\fullref{sec5}, with the added
simplification that we do not need to worry about the analytic quality of the
construction.

We let $D_{i,j}$ be a family of disjoint balls in $\SS^n$ for $i,j \in \Z$
such that there are homeomorphisms $T,F$ for which
$T$ takes $D_{i,j}$ to $D_{i+1,j}$ for all $i,j$, and $F$ takes
$D_{0,j}$ to $D_{0,j+1}$, and is the identity on $D_{i,j}$ when $i\ne 0$.

Let $Z_1,Z_2$ be homeomorphisms taking $D_{0,0}$ to $E_1$ and $E_2$ respectively.

For each $\ell \in \lbrace 1,\dots,6 \rbrace$ we define $f_\ell$ with support contained in 
the closure of the union of the $D_{i,j}$ by the formula
$$f_\ell = \prod_{i=0}^\infty \prod_{j=0}^\infty H_{i,\ell}^{Z_kF^{-j}T^{-i}}$$
where $k=1$ if $H_{i,\ell}$ has support in $E_1$, and $k=2$ if $H_{i,\ell}$ has support
in $E_2$.

Then as before, we can write $H_{i,\ell}$ as a word of length $\sim 4i$ in $f_1,\dots,f_6,T,F,Z_1,Z_2$
and their inverses. In detail:
$$H_{i,\ell} = Z_k y_{i,\ell} Z_k^{-1}$$
where
$$y_{i,\ell} = w_{i,\ell} F w_{i,\ell}^{-1} F^{-1}$$
and
$$w_{i,\ell} = T^{-i} f_\ell T^i.$$
Since we can do this for each $i,\ell$, we can exhibit each $h_i$
as a distortion element, whose distortion function grows as fast as desired.
Note that by choosing the $n_i$ to all be mutually coprime, we can ensure that
the $h_i$ are all actually contained in the group in question.
\end{proof}

\newpage

\appendix\def\theshortauthors{Yves de Cornulier}\def\theshorttitle{Appendix}
\section{Strong boundedness of $\homeo(S^n)$}
\cl{\sc Yves de Cornulier}

\bigskip
\begin{defn}
A group $G$ is strongly bounded\footnote{The following
terminologies for the same concept also exist in the literature:
$G$ has the Bergman Property; $G$ has the strong Bergman Property;
$G$ has uncountable strong cofinality.} if it satisfies one of the
following equivalent conditions:
\begin{itemize}
    \item[(i)] Every length function on $G$, ie function $L\co G\to\R_+$
        satisfying $L(1)=0$, $L(g^{-1})=L(g)$ and $L(gh)\le L(g)+L(h)$ for all $g,h\in G$, is bounded.

    \item[(ii)] Every action of $G$ by isometries on a metric space has
        bounded orbits.\vspace{10pt}

    \item[(iii)] \begin{itemize}
        \item[$\bullet$] $G$ is Cayley bounded: for every symmetric generating subset
            $S$ of $G$, there exists $n$ such that $G\subset
            S^n=\{s_1\dots s_n\,|\;s_1,\dots,s_n\in S\}$, and
        \item[$\bullet$] $G$ has uncountable cofinality, ie $G$ cannot be
            expressed as the union of an increasing sequence of proper
            subgroups.
        \end{itemize}
\end{itemize}\label{def:SB}
\end{defn}

The definition of groups with uncountable cofinality appeared in
the characterization by Serre \cite[Section~6.1]{Serre} of groups with
Property (FA), meaning that every isometric action on a simplicial
tree has a fixed point. For instance, a countable group has
uncountable cofinality if and only if it is finitely generated.

Much later, the concept of strong boundedness was introduced by
Bergman \cite{Bergman}, where it is proved that the permutation
group of any set is strongly bounded. Subsequently, intensive
research on the subject has been carried on (see, among others,
de Cornulier\cite{Cornulier}, Droste--Holland \cite{Droste_Holland},
Kechris--Rosendal \cite{Kechris_Rosendal} and Khelif \cite{Khelif}, and the
references in Bergman \cite{Bergman}). It is worth noting that a countable
group is strongly bounded if and only if it is finite, so that the
definition is of interest only for uncountable groups. In
\fullref{def:SB}, the equivalence between (i) and (ii) is
easy and standard; the equivalence between Conditions (i) and
(iii) is established in \cite{Cornulier} but already apparent
in~\cite{Bergman}.

Fix an integer $n\ge 1$. The purpose of this appendix is to point
out the following consequence of the proof of \fullref{thmc} in the
paper above.

\begin{thm}The group $\homeo(\SS^n)$ is strongly
bounded.\label{thm:homeo_Sn_SB}
\end{thm}
A weaker version of \fullref{thm:homeo_Sn_SB} was recently
proved by Rosendal \cite[Theorem 1.7]{Rosendal}; namely, $\homeo(\SS^n)$ is
strongly bounded as a \emph{topological} group (for the uniform
convergence); this means that every \emph{continuous} length
function is bounded. In \cite[Theorem 5.4]{Rosendal}, it was also
proved that $\homeo(\SS^1)$ is strongly bounded.
In contrast, if $M$ is any compact non-discrete manifold, then for every 
$r\ge 1$, the groups $\textnormal{Diff}^r(M)$ and their connected 
components are not strongly bounded. Indeed, the length function $h\mapsto 
\|dh\|^+$ defined in \fullref{general_nonsense} is unbounded.

\begin{proof}
Clearly, it suffices to show that the subgroup of index two
$\homeo^+(\SS^n)$ is strongly bounded. By contradiction, we suppose
the existence of an unbounded length function $L$ on $G$. Let us
pick a sequence $(h_i)$ in $G$ satisfying $L(h_i)\ge i^2$ for
all~$i$.

Using the notation in the proof of \fullref{thmc}, set
$S=\penalty0\{f_1,\dots,\penalty0 f_6,\penalty0 T,\penalty0 F,
Z_1,\penalty0 Z_2\}$. Then each $h_i$ can be expressed by a word
of length $\sim 24i$ in $S^\pm$. But, on the subgroup of
$\homeo^+(\SS^n)$ generated by the finite set $S$, the length
function $L$ must be dominated by the word length with respect to
$S$. This contradicts the assumption $L(h_i)\ge i^2$ for all~$i$.
\end{proof}

\begin{rmk}
A similar argument to that of the proof of \fullref{thm:homeo_Sn_SB} was used in \cite{Cornulier} to prove that
$\omega_1$--existentially closed groups are strongly bounded. This
reasoning is made systematic by Khelif \cite{Khelif}. Let us say that a group is
\textit{strongly distorted} (introduced as ``Property~P" in \cite{Khelif}) if
there exists an integer $m$ and an integer-valued sequence $(w_n)$
with the following property: for every sequence $(h_n)$ in $G$,
there exist $g_1,\dots,g_m\in G$ such that, for every $n$, one can
express $h_n$ as an element of length $w_n$ in the $g_i$'s.
Following the proof of \fullref{thm:homeo_Sn_SB}, we get that
a strongly distorted group is strongly bounded, and that
$\homeo(\SS^n)$ is strongly distorted.

The symmetric group on any set is strongly distorted: a proof can
be found in Galvin \cite{Galvin}, although a weaker result is stated
there. In \cite{Khelif}, it is claimed that the automorphism group
of any 2--transitive chain is strongly distorted; strong
boundedness was previously proved by Droste and Holland
\cite{Droste_Holland}. On the other hand, if $F$ is a nontrivial finite
perfect group, then the infinite (unrestricted) direct product $F^\N$
is strongly bounded \cite{Cornulier}; however it is clearly not strongly
distorted since it is infinite and locally finite.
\end{rmk}

\begin{rmk}
It follows from \fullref{thm:homeo_Sn_SB} that $\homeo(\SS^n)$
is Cayley bounded (see \fullref{def:SB}); that is, the
Cayley graph with respect to any generating subset has bounded
diameter. It is natural to ask whether there is a uniform bound on
those diameters: the answer is negative. Indeed, endow $\SS^n$ with
its Euclidean metric, and, for $r>0$, set
$$W_r=\{g\in\homeo(\SS^n)\,|\;\forall x\in \SS^n,\,d(x,g(x))<r\}$$
and
$$W_r^+=W_r\cap\homeo^+(\SS^n).$$
Then $W_r^+$ is open in
$\homeo^+(\SS^n)$. The group $\homeo^+(\SS^n)$, endowed with the
topology of uniform convergence, is connected: this is well-known,
and can be deduced, for instance, from \fullref{KS_lemma} above.
It follows that $W_r^+$ generates $\homeo^+(\SS^n)$. Clearly, for
every $k\ge 1$, we have $(W_r^+)^k\subset W_{kr}^+$. It follows
that if we have chosen $r\le 2/k$, then
$(W_r^+)^k\neq\allowbreak\homeo^+(\SS^n)$. Thus $\homeo^+(\SS^n)$ has
Cayley graphs of arbitrary large diameter. A similar argument
works for $\homeo(\SS^n)$ as follows: fix a reflection
$T\in\textnormal{O}(n+\penalty10000 1)$ of $\SS^n$, and take
$W'_r=W_r\cup TW_r$. Then it is easy to check that
$(W'_r)^k\subset W'_{kr}$ for all $k\ge 1$, so that
$(W'_r)^k\neq\homeo(\SS^n)$ if we have chosen $r\le 2/k$.
\end{rmk}

\def\theshorttitle{References}
\bibliographystyle{gtart}
\bibliography{link}

\end{document}